\documentclass[12pt]{article}

\usepackage{amsmath,amssymb,amsthm,amscd}
\usepackage[colorlinks, linkcolor=blue, anchorcolor=gree, citecolor=red]{hyperref}
\usepackage[numbers,sort&compress]{natbib}
\usepackage{graphicx}
\usepackage{amssymb}
\usepackage{mathrsfs}
\usepackage{amssymb, amsmath, fancyhdr, graphicx, graphics, pstricks}

\begin{document}

\newcommand{\Cyc}{{\rm{Cyc}}}\newcommand{\diam}{{\rm{diam}}}

\newtheorem{thm}{Theorem}[section]
\newtheorem{pro}[thm]{Proposition}
\newtheorem{lem}[thm]{Lemma}
\newtheorem{con}[thm]{Conjecture}
\newtheorem{fac}[thm]{Fact}
\newtheorem{cor}[thm]{Corollary}
\theoremstyle{definition}
\newtheorem{ex}[thm]{Example}
\newtheorem{ob}[thm]{Observtion}
\newtheorem{remark}[thm]{Remark}
\newcommand{\bth}{\begin{thm}}
\renewcommand{\eth}{\end{thm}}
\newcommand{\bex}{\begin{examp}}
\newcommand{\eex}{\end{examp}}
\newcommand{\bre}{\begin{remark}}
\newcommand{\ere}{\end{remark}}

\newcommand{\bal}{\begin{aligned}}
\newcommand{\eal}{\end{aligned}}
\newcommand{\beq}{\begin{equation}}
\newcommand{\eeq}{\end{equation}}
\newcommand{\ben}{\begin{equation*}}
\newcommand{\een}{\end{equation*}}

\newcommand{\bpf}{\begin{proof}}
\newcommand{\epf}{\end{proof}}
\renewcommand{\thefootnote}{}
\newcommand{\sdim}{{\rm sdim}}

\def\beql#1{\begin{equation}\label{#1}}
\title{\Large\bf Two bounds for generalized $3$-connectivity of Cartesian product graphs}

\author{{\sc Hui Gao~~~~~~ Benjian Lv ~~~~~~Kaishun Wang
}\\[15pt]
{\small\em }Sch. Math. Sci. {\rm \&} Lab. Math. Com. Sys., \\
Beijing Normal University, Beijing, 100875, China
\\
}

 \date{}

\maketitle

\begin{abstract}
The generalized $k$-connectivity $\kappa_{k}(G)$ of a graph $G$, which was introduced by  Chartrand et al.(1984) is a generalization of the concept of vertex connectivity. Let $G$ and $H$ be nontrivial connected graphs. Recently, Li et al. gave a lower bound for the generalized $3$-connectivity of the Cartesian product graph $G \square H$ and proposed a conjecture for the case that $H$ is $3$-connected. In this paper, we give two different forms of lower bounds for the generalized $3$-connectivity of Cartesian product graphs. The first lower bound is stronger than theirs, and the second confirms their conjecture.
\end{abstract}



\footnote{E-mail addresses: gaoh1118@yeah.net (H. Gao), bjlv@bnu.edu.cn (B. Lv), \\wangks@bnu.edu.cn (K. Wang).}
\section{Introduction}
All graphs in this paper are undirected, finite and simple. We refer to the book \cite{bondy} for graph theoretic notations and terminology not described here. The generalized connectivity of a graph G, which was introduced by Chartrand et al. \cite{chartrand}, is a natural generalization of the concept of vertex connectivity.

A tree $T$ is called an $S$-tree if $S \subseteq V(T)$. A family of $S$-trees $T_{1}$, $T_{2}$,...,$T_{r}$ are internally disjoint if $E(T_{i}) \cap E(T_{j}) = \phi$ and $V(T_{i}) \cap V(T_{j}) = S$ for any pair of integers $i$ and $j$, where $1 \leq i < j \leq r$. We denote by $\kappa(S)$ the greatest number of internally disjoint $S$-trees. For an integer $k$ with $2 \leq k \leq v(G)$, the generalized $k$-connectivity $\kappa_{k}(G)$ are defined to be the least value of $\kappa (S)$ when $S$ runs over all $k$-subsets of $V(G)$. Clearly, when $k=2$, $\kappa_{2}(G)=\kappa(G)$.

In addition to being a natural combinatorial notation, the generalized connectivity can be motivated by its interesting interpretation in practice. For example, suppose that $G$ represents a network. If one considers to connect a pair of vertices of $G$, then a path is used to connect them. However, if one wants to connect a set $S$ of vertices of $G$ with $|S| \geq 3$, then a tree has to be used to connect them. This kind of tree with minimum order for connecting a set of vertices is usually called Steiner tree, and popularly used in the physical design of VLSI, see \cite{NA}. Usually, one wants to consider how tough a network can be, for the connection of a set of vertices. Then, the number of totally independent ways to connect them is a measure for this purpose. The generalized $k$-connectivity can serve for measuring the capability of a network $G$ to connect any $k$ vertices in $G$.

Determining  $\kappa_{k}(G)$ for most graphs is a difficult problem. In \cite{shasha}, Li and Li derived that for any fixed integer $k \geq 2$, given a graph $G$ and a subset $S$ of $V(G)$, deciding whether there are $k$ internally disjoint trees connecting $S$, namely deciding whether $\kappa(S) \geq k$ is NP-complete. The exact value of $\kappa_{k}(G)$ is known for only a small class of graphs. Examples are complete graphs \cite{complete}, complete bipartite graphs \cite{complete bipartite}, complete equipartition $3$-partite graphs \cite{complete equipartition 3-partite}, star graphs \cite{star}, bubble-sort graphs \cite{star}, and connected Cayley graphs on Abelian groups with small degrees \cite{small cayley}. Upper bounds and lower bounds of generalized connectivity of a graph have been investigated by Li et al. \cite{hengzhe, bound connectivity, bound edge-connectivity} and Li and Mao \cite{lexico}. And Li et al. investigated Extremal problems in \cite{minimal, minimal size}. We refer the readers to \cite{survey} for more results.

In \cite{hengzhe}, Li et al studied the generalized $3$-connectivity of Cartesian product graphs and showed the following result.

\begin{thm} [\cite{hengzhe}] \label{1}
Let $G$ and $H$ be connected graphs such that $\kappa_{3}(G) \geq \kappa_{3}(H)$. The following assertions hold:

(i) if $\kappa(G)=\kappa_{3}(G)$, then $\kappa_{3}(G \square H) \geq \kappa_{3}(G) + \kappa_{3}(H) -1$. Moreover, the bound is sharp;

(ii) if $\kappa(G) > \kappa_{3}(G)$, then $\kappa_{3}(G \square H) \geq \kappa_{3}(G) + \kappa_{3}(H)$. Moreover, the bound is sharp.
\end{thm}

Later in \cite{graph products}, Li et al gave a better result when $H$ becomes a $2$-connected graph.

\begin{thm} [\cite{graph products}]  \label{2}
Let $G$ be a nontrivial connected graph, and let $H$ be a 2-connected graph. The following assertions hold:

(i) if $\kappa(G)=\kappa_{3}(G)$, then $\kappa_{3}(G \square H) \geq \kappa_{3}(G) + 1$. Moreover, the bound is sharp;

(ii) if $\kappa(G) > \kappa_{3}(G)$, then $\kappa_{3}(G \square H) \geq \kappa_{3}(G) + 2$. Moreover, the bound is sharp.
\end{thm}

Also in \cite{graph products}, Li et al proposed a conjecture as follows:

\begin{con} [\cite{graph products}] \label{3}
Let $G$ be a nontrivial connected graph, and let $H$ be a $3$-connected graph. The following assertions hold:

(i) if $\kappa(G)=\kappa_{3}(G)$, then $\kappa_{3}(G \square H) \geq \kappa_{3}(G) + 2$. Moreover, the bound is sharp;

(ii) if $\kappa(G) > \kappa_{3}(G)$, then $\kappa_{3}(G \square H) \geq \kappa_{3}(G) + 3$. Moreover, the bound is sharp.
\end{con}

In this paper, we give two different forms of lower bounds for generalized $3$-connectivity of Cartesian product graphs.

\begin{thm}\label{15}
Let $G$ and $H$ be nontrivial connected graphs. Then $\kappa_{3}(G \square H) \geq min \{ \kappa_{3}(G)+ \delta(H)$, $\kappa_{3}(H)+ \delta(G)$, $\kappa(G) + \kappa(H)-1\}$.
\end{thm}

\begin{thm}\label{16}
Let $G$ be a nontrivial connected graph, and let $H$ be an $l$-connected graph. The following assertions hold:

(i) if $\kappa(G)=\kappa_{3}(G)$ and $1 \leq l \leq 7$, then $\kappa_{3}(G \square H) \geq \kappa_{3}(G) + l -1$. Moreover, the bound is sharp;

(ii) if $\kappa(G) > \kappa_{3}(G)$ and $1 \leq l \leq 9$, then $\kappa_{3}(G \square H) \geq \kappa_{3}(G) + l$. Moreover, the bound is sharp.
\end{thm}

The paper is organized as follows. In Section $2$, we introduce some definitions and notations. In Section $3$, we give a proof of theorem~\ref{15}, which induces theorem~\ref{1} and theorem~\ref{2}, and confirms conjecture~\ref{3}. In section $4$, we discuss the problem which number the connectivity of $H$ can be such that conjecture~\ref{3} still holds. And theorem~\ref{16} is our answer and there are counterexamples when $l \geq 8$ for $\kappa(G) = \kappa_{3}(G)$ and $l \geq 10$ for $\kappa(G) > \kappa_{3}(G)$.
\section{Preliminaries}
Let $G$ and $H$ be two graphs with $V(G)= \{ u_{1}, u_{2},..., u_{n} \}$ and $V(H)= \{ v_{1}, v_{2},..., v_{m} \}$, respectively. Let $\kappa(G)=k$, $\kappa(H)=l$, $\delta(G)=\delta_{1}$, and $\delta(H)=\delta_{2}$. And the discussion below is always based on the hypotheses.

Recall that the Cartesian product (also called the square product) of two graphs $G$ and $H$, written as $G \square H$,  is the graph with vertex set $V(G) \times V(H)$, in which two vertices $( u, v )$ and $(u', v')$ are adjacent if and only if $u = u'$ and $vv' \in E(H)$, or $v=v'$ and $uu' \in E(G)$. By starting with a disjoint union of two graphs $G$ and $H$ and adding edges joining every vertex of $G$ to every vertex of $H$, one obtains the join of $G$ and $H$, denoted by $G \vee H$.

For any subgraph $G_{1} \subseteq G $, we use $G_{1}^{v_{j}}$ to denote the subgraph of $G \square H$ with vertex set $\{ (u_{i}, v_{j}) |u_{i} \in V(G_{1}) \}$ and edge set $\{ (u_{i_{1}}, v_{j})(u_{i_{2}}, v_{j}) |u_{i_{1}} u_{i_{2}} \in E(G_{1}) \}$. Similarly, for any subgraph $H_{1} \subseteq H$, we use $H_{1}^{u_{i}}$ to denote the subgraph of $G \square H$ with vertex set $\{ (u_{i}, v_{j}) |v_{j} \in V(H_{1}) \}$ and edge set $\{ (u_{i}, v_{j_{1}})(u_{i}, v_{j_{2}}) |v_{j_{1}} v_{j_{2}} \in E(H_{1}) \}$. Clearly, $G_{1}^{v_{j}} \cong G_{1}$, $H_{1}^{u_{i}} \cong H_{1}$.

Let $x \in V(G)$ and $Y \subseteq V(G)$. An $(x$, $Y)$-path is a path which starts at x, ends at a vertex of $Y$, and whose internal vertices do not belong to $Y$. A family of $k$ internally disjoint $(x$, $Y)$-paths whose terminal vertices are distinct is referred to as a $k$-fan from $x$ to $Y$.

For some $1 \leq t \leq \lfloor \frac{k}{2} \rfloor$ and $s \geq t+1$, in $G$, a family $\{P_{1}, P_{2}, ..., P_{s} \}$ of $s$ $u_{1}u_{2}-$paths is called an $(s, t)$-original-path-bundle with respect to $(u_{1}, u_{2}, u_{3})$, if $u_{3}$ are on $t$ paths $P_{1}$, ..., $P_{t}$, and the $s$ paths have no internal vertices in common except $u_{3}$, as shown in figure~\ref{f8}.a. If there is not only an $(s,t)$-original-path-bundle $\{P_{1}^{'}, P_{2}^{'}, ..., P_{s}^{'} \}$ with respect $(u_{1}, u_{2}, u_{3})$, but also a family $\{M_{1}, M_{2}, ..., M_{k-2t} \}$ of $k-2t$ internally disjoint $(u_{3}, X)$-paths avoiding the vertices in $V(P_{1}^{'} \cup ... \cup P_{t}^{'}) - \{ u_{1}, u_{2}, u_{3}\}$, where $X = V(P_{t+1}^{'} \cup ... \cup P_{s}^{'})$, then we call the family of paths $\{ P_{1}^{'}, P_{2}^{'}, ..., P_{s}^{'} \} \cup \{M_{1}, M_{2}, ..., M_{k-2t} \}$ an $(s, t)$-reduced-path-bundle with respect to $(u_{1}, u_{2}, u_{3})$, as shown in figure~\ref{f8}.b.

\begin{figure}[hptb]
  \centering
  \includegraphics[width=17cm]{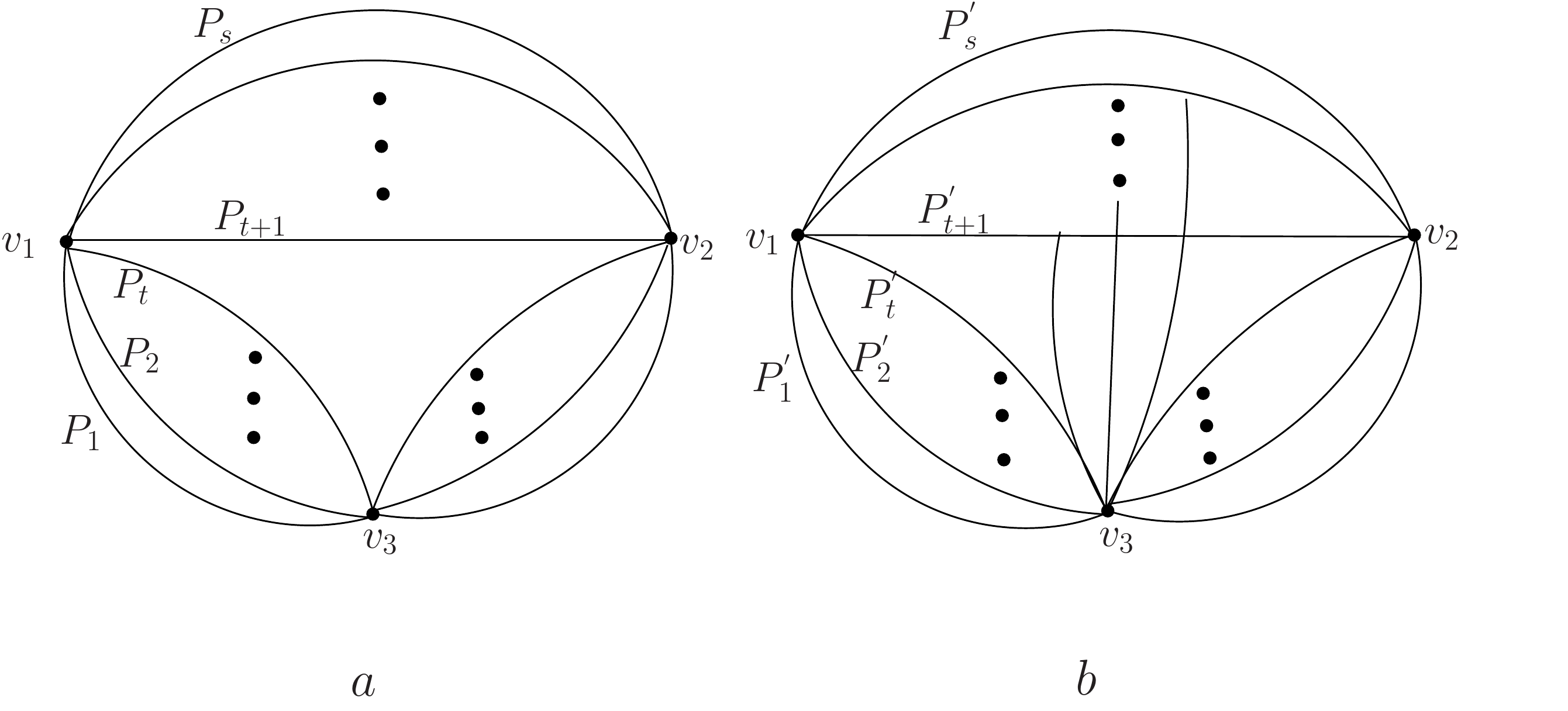}\\
  \caption{}\label{f8}
\end{figure}

In order to show our main results, we need the following theorems and lemmas.

\begin{lem}{\rm(\cite[Fan Lemma]{bondy})} \label{4}
Let $G$ be a $k$-connected graph, let $x$ be a vertex of $G$, and let $Y \subseteq V(G)\setminus \{ x\}$. Then there exists a $k$-fan in $G$ from $x$ to $Y$.
\end{lem}

\begin{thm}{\rm(\cite[p.219]{bondy})} \label{5}
Let $S$ be a set of three pairwise-nonadjacent edges in a simple $3$-connected graph $G$. Then there is a cycle in $G$ containing all three edges of $S$ unless $S$ is an edge cut of $G$.
\end{thm}

\begin{thm}[\cite{bound connectivity}] \label{6}
Let $G$ be a connected graph with at least three vertices. If $G$ has two adjacent vertices with minimum degree $\delta$, then $\kappa_{3}(G) \leq \delta -1$.
\end{thm}

\begin{thm}[\cite{bound connectivity}] \label{7}
Let $G$ be $k$-connected, and $u_{1}, u_{2}, u_{3} \in V(G)$. Then for some $0 \leq t \leq \lfloor \frac{k}{2} \rfloor $, there exists a $(k, t)$-reduced-path-bundle $\{ P_{1}, P_{2}, ..., P_{k}\} \cup \{ M_{1}, M_{2}, ..., M_{k-2t} \}$ such that for any $1 \leq i \leq k-2t$, the terminal vertex of $M_{i}$ is on $P_{t+i}$.
\end{thm}

\begin{thm}[\cite{bound connectivity}] \label{8}
let $G$ be a connected graph with $n$ vertices. For every two integers $k$ and $r$ with $k \geq 0$ and $r \in \{ 0$, $1$, $2$, $3 \}$, if $\kappa(G) = 4k +r$, then $3k + \lceil \frac{r}{2}\rceil \leq \kappa_{3}(G) \leq \kappa (G) $. Moreover, the bound is sharp.
\end{thm}

\begin{lem}[\cite{complete bipartite}]\label{9}
Let $a$, $b$ be integers such that $a+b \geq 3$ and $a \leq b$. Then,

\begin{equation}
\kappa_{3}(K_{a,b}) =
   \begin{cases}
   a-1& \text{if $a=b$},\\
   a& \text{if $a < b$}.
   \end{cases}
\end{equation}

\end{lem}

\begin{lem}[\cite{graph products}]\label{10}
Let $C_{1}$, $C_{2}$, ..., $C_{k}$ be cycles. then $\kappa_{3}(C_{1} \square C_{2} \square ... \square C_{k}) = 2k-1$.
\end{lem}

\section{One lower bound}

\begin{lem}\label{11}
If $S= \{ (u_{1}$, $v_{1})$, $(u_{2}$, $v_{2})$, $(u_{3}$, $v_{3})\}$, then $\kappa(S) \geq k+l-1$.
\end{lem}

\proof
Set $S^{'} = \{ (u_{i}$, $v_{j}) | i$, $j = 1$, $2$, $3\}$.

Case 1 $(G \square H) (S^{'})$ is not isomorphic to $ C_{3} \square C_{3}$.

Since $(G \square H) (S^{'})$ is not isomorphic to $ C_{3} \square C_{3}$, either $G ( \{ u_{1}$, $u_{2}$, $u_{3} \} )$ or $H ( \{ v_{1}$, $v_{2}$, $v_{3} \} )$ is not isomorphic $C_{3}$. Without loss of generality, suppose $v_{1}$ is not adjacent to $v_{2}$ in $H$. Since $H$ is $l$-connected, there exist $l$ internally disjoint paths in H from $v_{1}$ to $v_{2}$, say $P_{j}$, in which $v_{i_{j}}$ is adjacent $v_{2}$,  $j=1$, $2$, ..., $l$. Suppose $v_{3}$ is not in $P_{j}$, $j=1$, $2$, ..., $l-1$. Also, according to Lemma~\ref{4}, there exists an $l$-fan in $H$ from $v_{3}$ to $\{ v_{i_{1}}$, $v_{i_{2}}$, ..., $v_{i_{l-1}}$, $v_{1} \}$, say $Q_{j}$, which is a $v_{3}v_{i_{j}}$-path, $j=1$, $2$, ..., $l-1$; and $Q_{l}$, which is a $v_{3}v_{1}$-path. Set $T_{j}=(P_{j}-v_{2})^{u_{1}} \cup G^{v_{i_{j}}} \cup (v_{2}v_{i_{j}})^{u_{2}} \cup Q_{j}^{u_{3}}$, $j=1$, $2$, ..., $l-1$. It is easy notice that, for $j=1$, $2$, ..., $l-1$, $(P_{j}-v_{2})^{u_{1}}$ connects $(u_{1},v_{1})$ and $(u_{1},v_{i_{j}})$; $G^{v_{i_{j}}} \cup (v_{2}v_{i_{j}})^{u_{2}}$ connects $(u_{1},v_{i_{j}})$, $(u_{2},v_{2})$ and $(u_{3},v_{i_{j}})$; and $Q_{j}^{u_{3}}$ connects $(u_{3},v_{i_{j}})$ and $(u_{3},v_{3})$. Thus, $T_{j}$ connects $(u_{1},v_{1})$, $(u_{2},v_{2})$ and $(u_{3},v_{3})$. Since $G$ is $k$-connected, there exist $k$ internally disjoint paths in $G$ from $u_{1}$ to $u_{2}$, say $R_{j}$, in which $u_{i_{j}^{'}}$ is adjacent to $u_{2}$, $j=1$, $2$, ..., $k$. Suppose $u_{i_{1}^{'}}$, $u_{i_{2}^{'}}$, ..., $u_{i_{k-2}^{'}} \neq u_{1} $ or $u_{3}; u_{i_{k-1}^{'}} \neq u_{3}$; and $u_{i_{k}^{'}} \neq u_{1}$.

Case 1.1 $u_{i_{k}^{'}} \neq u_{3}$

Since $G$ is $k$-connected, there exists a $k$-fan in $G$ from $u_{3}$ to $\{ u_{i_{1}^{'}}$, $u_{i_{2}^{'}}$, ..., $u_{i_{k-2}^{'}}$, $u_{2}$, $u_{i_{k}^{'}} \}$, say $S_{j}$, which is a $u_{3}u_{i_{j}^{'}}$-path, $j=1$, $2$, ..., $k-2$, $k$; and $S_{k-1}$, which is a $u_{3}u_{2}$-path. Since $H$ is $l$-connected, $H - \{ v_{i_{1}}$, $v_{i_{2}}$, ..., $v_{i_{l-1}} \}$ is connected. Set $T_{j}^{'}=(R_{j}-u_{2})^{v_{1}} \cup (H - \{ v_{i_{1}}, v_{i_{2}}, ...,  v_{i_{l-1}} \})^{ u_{i_{j}^{'}}} \cup (u_{i_{j}^{'}}u_{2})^{v_{2}} \cup S_{j}^{v_{3}}$, $j=1$, $2$, ..., $k-2$, $k$; and $T_{k-1}^{'}= R_{k-1}^{v_{1}} \cup (H - \{ v_{i_{1}}, v_{i_{2}}, ...,  v_{i_{l-1}} \})^{ u_{2}} \cup (S_{k-1})^{v_{3}}$. It is easy to notice that, for $j=1$, $2$, ..., $k-2$, $k$, $(R_{j}-u_{2})^{v_{1}}$ connects $(u_{1}$, $v_{1})$ and $(u_{i_{j}^{'}}$, $v_{1})$; $(H - \{ v_{i_{1}}, v_{i_{2}}, ...,  v_{i_{l-1}} \})^{ u_{i_{j}^{'}}} \cup (u_{i_{j}^{'}}u_{2})^{v_{2}}$ connects $(u_{i_{j}^{'}}$, $v_{1})$, $(u_{2}$, $v_{2})$ and $(u_{i_{j}^{'}}$, $v_{3})$; and $S_{j}^{v_{3}}$ connects $((u_{i_{j}^{'}}$, $v_{3}))$ and $(u_{3}$, $v_{3})$. Thus $T_{j}^{'}$ connects $(u_{1}$, $v_{1})$, $(u_{2}$, $v_{2})$ and $(u_{3}$, $v_{3})$. Similarly, so does $T_{k-1}^{'}$. It is clear that $T_{1}$, ..., $T_{l-1}$, $T_{1}^{'}$, ..., $T_{k}^{'}$ are pairwise disjoint except the vertex set $S$. See figure~\ref{f1}.

\begin{figure}[hptb]
  \centering
  \includegraphics[width=10cm]{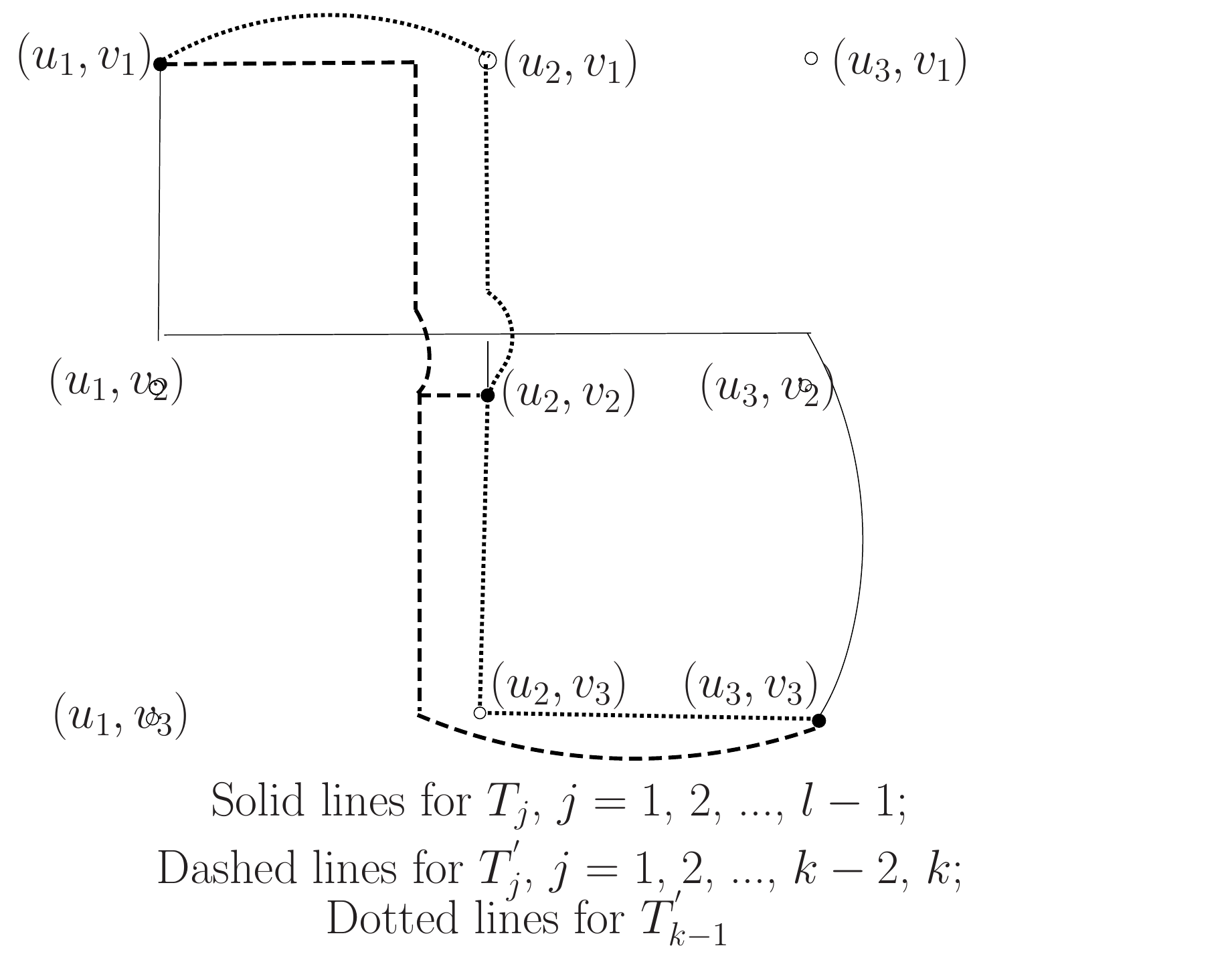}\\
  \caption{}\label{f1}
\end{figure}

Case 1.2 $u_{i_{k}^{'}} = u_{3}$

Since $G$ is $k$-connected, there exists a $(k-1)$-fan in $G-u_{1}$ from $u_{3}$ to $\{ u_{i_{1}^{'}}$, $u_{i_{2}^{'}}$, ..., $u_{i_{k-2}^{'}}$, $u_{2} \}$, say $S_{j}$, which is a $u_{3}u_{i_{j}^{'}}$-path, $j=1$, $2$, ..., $k-2$; and $S_{k-1}$, which is a $u_{3}u_{2}$-path. Set $T_{j}^{'}=(R_{j}-u_{2})^{v_{1}} \cup (H - \{ v_{i_{1}}, v_{i_{2}}, ...,  v_{i_{l-1}} \})^{ u_{i_{j}^{'}}} \cup (u_{i_{j}^{'}}u_{2})^{v_{2}} \cup S_{j}^{v_{3}}$, $j=1$, $2$, ..., $k-2$; $T_{k-1}^{'}= R_{k-1}^{v_{1}} \cup (H - \{ v_{i_{1}}$, $v_{i_{2}}$, ...,  $v_{i_{l-1}} \})^{ u_{2}} \cup S_{k-1}^{v_{3}}$; and $T_{k}^{'}=Q_{l}^{u_{3}} \cup (R_{k}-u_{2})^{v_{1}} \cup P_{l}^{u_{1}} \cup R_{k-1}^{v_{2}}$. It is easy to notice that $Q_{l}^{u_{3}}$ connects $(u_{3}$, $v_{3})$ and $(u_{3}$, $v_{1})$;  $(R_{k}-u_{2})^{v_{1}}$ connects $(u_{3}$, $v_{1})$ and $(u_{1}$, $v_{1})$; $P_{l}^{u_{1}}$ connects $(u_{1}$, $v_{1})$ and $(u_{1}$, $v_{2})$; and $R_{k-1}^{v_{2}}$ connects $(u_{1}$, $v_{2})$ and $(u_{2}$, $v_{2})$. Thus, $T_{k}^{'}$ connects $(u_{1}$, $v_{1})$, $(u_{2}$, $v_{2})$ and $(u_{3}$, $v_{3})$. It is clear that $T_{1}$, ..., $T_{l-1}$, $T_{1}^{'}$, ..., $T_{k}^{'}$ are pairwise disjoint except the vertex set $S$. See figure~\ref{f2}.

\begin{figure}[hptb]
  \centering
  \includegraphics[width=10cm]{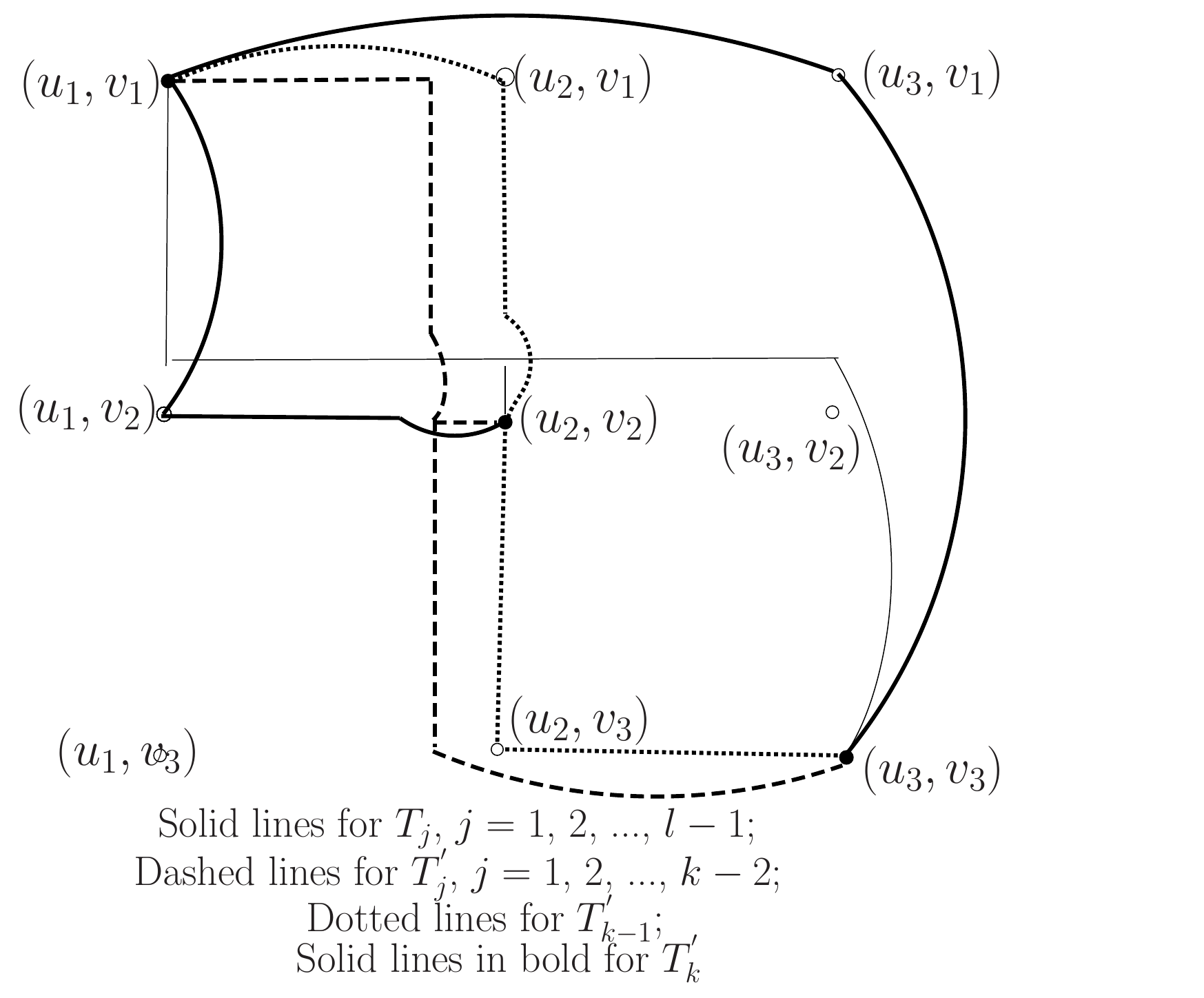}\\
  \caption{}\label{f2}
\end{figure}

Case 2 $(G \square H) (S^{'}) \cong C_{3} \square C_{3}$.

Because of lemma~\ref{10}, there exist $3$ internally disjoint $S$-trees in $(G \square H) (S^{'})$, say $T_{j}^{''}$, $j=1$, $2$, $3$. Since $H$ is $l$-connected, $\kappa$ $((H - v_{3})- v_{1}v_{2}) \geq l-2$. Thus there exist $l-2$ internally disjoint $v_{1}v_{2}$-paths in $(H - v_{3})- v_{1}v_{2}$, say $P_{j}$, in which $v_{i_{j}}$ is adjacent of $v_{2}$, $j=1$, $2$, ..., $l-2$. By Lemma~\ref{4}, there exists an $(l-2)$-fan in $H-v_{1}-v_{2}$ from $v_{3}$ to $\{ v_{i_{1}}$, $v_{i_{2}}$, ..., $v_{i_{l-2}} \}$. Similarly to Case 1, we can construct $l-2$ internally disjoint $S$-trees $T_{1}$, $T_{2}$, ..., $T_{l-2}$. Also, because $\kappa (G)=k$, another $k-2$ $S$-trees $T_{1}^{'}$, $T_{2}^{'}$, ..., $T_{k-2}^{'}$ can be constructed. It is clear that $T_{1}$, ..., $T_{l-2}$, $T_{1}^{'}$, ..., $T_{k-2}^{'}$, $T_{1}^{''}$, $T_{2}^{''}$, $T_{3}^{''}$ are internally disjoint $S$-trees.
\qed

\medskip

\begin{lem}\label{12}
If $S= \{ (u_{1}$, $v_{1})$, $(u_{1}$, $v_{2})$, $(u_{2}$, $v_{1})\}$, then $\kappa(S) \geq k+l-1$.
\end{lem}

\proof
Since $H$ is $l$-connected, there exist $l$ internally disjoint $v_{1}v_{2}$-paths in $H$, say $P_{j}$, in which $v_{i_{j}}$ is adjacent to $v_{1}$, $j=1$, $2$, ..., $l$. Suppose $v_{i_{j}} \neq v_{2}$, $j=1$, $2$, ..., $l-1$. Set $T_{j}=P_{j}^{u_{1}} \cup G^{v_{i_{j}}} \cup (v_{1}v_{i_{j}})^{u_{2}}$, $j=1$, $2$, ..., $l-1$. Since $G$ is $k$-connected, there exist $k$ internally disjoint $u_{1}u_{2}$-paths in $G$, say $Q_{j}$, in which $u_{i_{j}^{'}}$ is adjacent to $u_{1}$, $j=1$, $2$, ..., $k$. Suppose $u_{i_{j}^{'}} \neq u_{2}$, $j=1$, $2$, ..., $k-1$. Set $T_{j}^{'}=Q_{j}^{v_{1}} \cup (H- \{v_{i_{1}}, v_{i_{2}}, ..., v_{i_{l-1}} \})^{u_{i_{j}^{'}}} \cup (u_{1}u_{i_{j}^{'}})^{v_{2}}$, $j=1$, $2$, ..., $k-1$; and $T^{''}=P_{l}^{u_{1}} \cup Q_{k}^{v_{1}}$. It is clear that $T_{1}$, ..., $T_{l-1}$, $T_{1}^{'}$, ..., $T_{k-1}^{'}$, $T^{''}$ are connected graphs containing $S$ and are pairwise disjoint except $S$. See figure~\ref{f3}. \qed

\medskip

\begin{figure}[hptb]
  \centering
  \includegraphics[width=9cm]{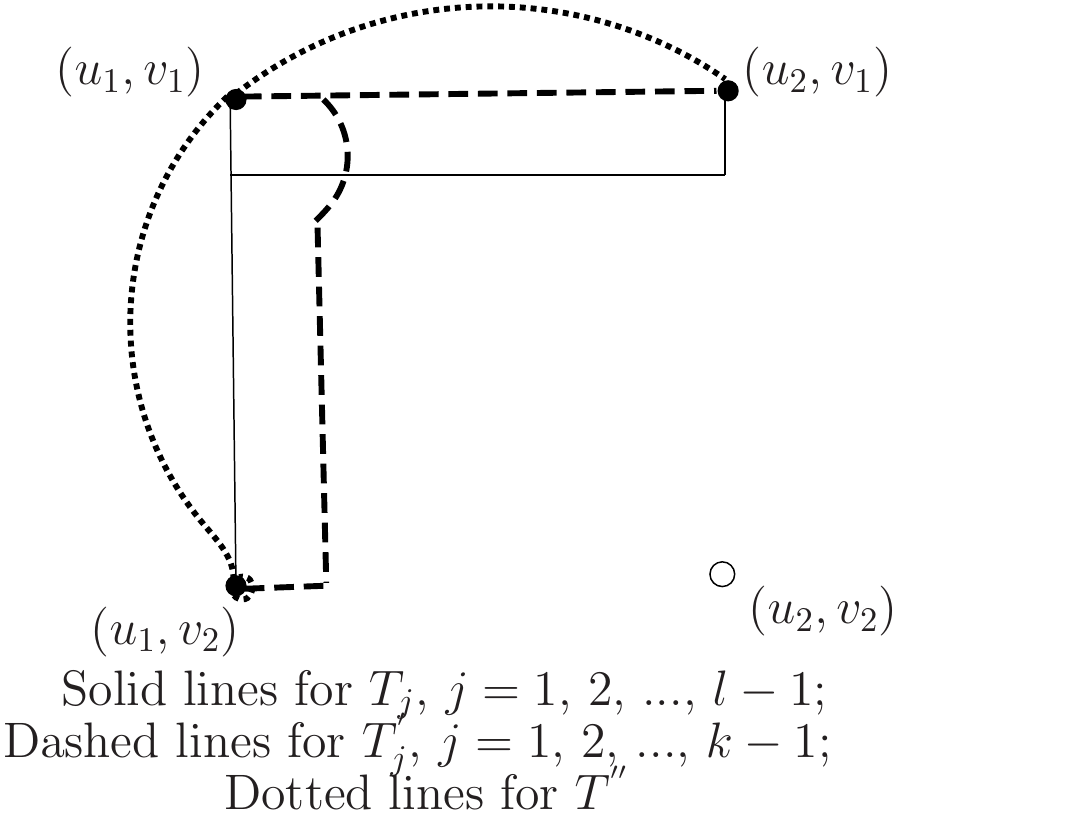}\\
  \caption{}\label{f3}
\end{figure}

\begin{lem}\label{13}
If $S= \{ (u_{1}$, $v_{1})$, $(u_{2}$, $v_{1})$, $(u_{3}$, $v_{2})\}$, then $\kappa (S) \geq k+l-1$.
\end{lem}

\proof
Since $H$ is $l$-connected, there exist $l$ internally disjoint $v_{1}v_{2}$-paths in $H$, say $P_{j}$, in which $v_{i_{j}}$ is adjacent to $v_{1}$, $j=1$, $2$, ..., $l$. Suppose $v_{i_{j}} \neq v_{2}$, $j=1$, $2$, ..., $l-1$. Set $T_{j}= (v_{1}v_{i_{j}})^{u_{1}} \cup (v_{1}v_{i_{j}})^{u_{2}} \cup G^{v_{i_{j}}} \cup (P_{j}-v_{1})^{u_{3}}$, $j=1$, $2$, ..., $l-1$. Since $G$ is $k$-connected, there exist $k$ internally disjoint $u_{1}u_{2}$-paths in $G$, say $Q_{j}$, in which $u_{i_{j}^{'}}$ is adjacent $u_{1}$, $j=1$, $2$, ..., $k$. Suppose $u_{i_{j}^{'}} \neq u_{2}$, $u_{3}$, $j=1$, $2$, ..., $k-2$. Due to Lemma~\ref{4}, there exists a $k$-fan in $G$ from $u_{3}$ to $\{ u_{i_{1}^{'}}$, $u_{i_{2}^{'}}$, ..., $u_{i_{k-2}^{'}}$, $u_{1}$, $u_{2}\}$, say $R_{j}$, which is a $u_{3}u_{i_{j}^{'}}$-path, $j=1$, $2$, ..., $k-2$; $R_{k-1}$, which is a $u_{3}u_{1}$-path; and $R_{k}$, which is a $u_{3}u_{2}$-path. Set $T_{j}^{'}=Q_{j}^{v_{1}} \cup (H - \{ v_{i_{1}}$, $v_{i_{2}}$, ..., $v_{i_{l-1}}\})^{u_{i_{j}^{'}}} \cup R_{j}^{v_{2}}$, $j=1$, $2$, ..., $k-2$; $T_{k-1}^{'}=Q_{k-1}^{v_{1}} \cup (H - \{ v_{i_{1}}$, $v_{i_{2}}$, ..., $v_{i_{l-1}}\})^{u_{1}} \cup R_{k-1}^{v_{2}}$; and $T_{k}^{'}=Q_{k}^{v_{1}} \cup (H - \{ v_{i_{1}}$, $v_{i_{2}}$, ..., $v_{i_{l-1}}\})^{u_{2}} \cup R_{k}^{v_{2}}$. It is clear that $T_{1}$, ..., $T_{l-1}$, $T_{1}^{'}$, ..., $T_{k}^{'}$ are connected graphs containing $S$ and are pairwise disjoint except $S$. See figure~\ref{f4}. \qed

\medskip

\begin{figure}[hptb]
  \centering
  \includegraphics[width=13cm]{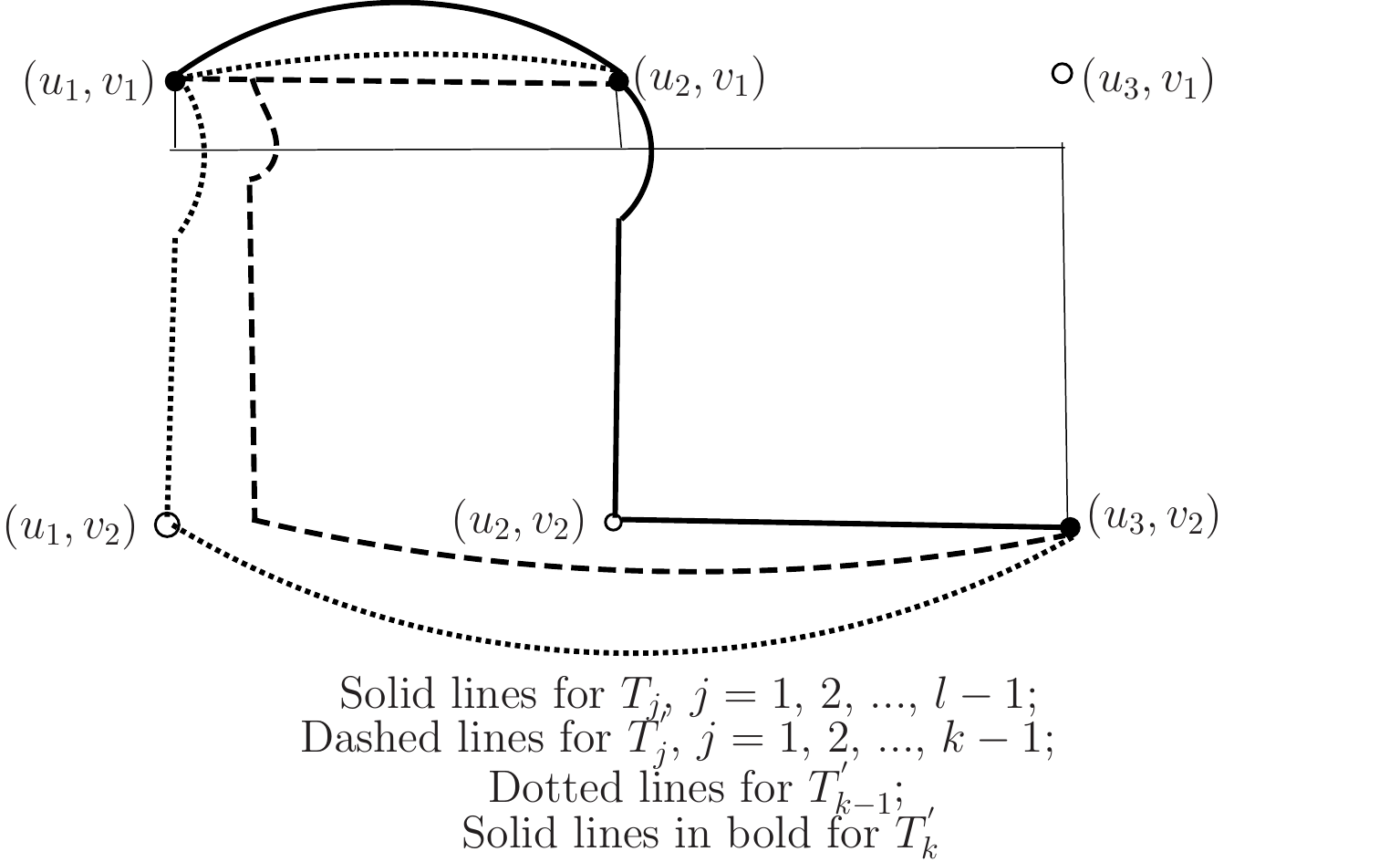}\\
  \caption{}\label{f4}
\end{figure}

\begin{lem}\label{14}
If $S= \{ (u_{1}$, $v_{1})$, $(u_{2}$, $v_{1})$, $(u_{3}$, $v_{1})\}$, then $\kappa (S) \geq \kappa_{3}(G)+ \delta(H)$.
\end{lem}

\proof
Let $T_{1}$, $T_{2}$, ..., $T_{\kappa_{3}(G)}$ be the internally disjoint $\{ u_{1}$, $u_{2}$, $u_{3}\}$-trees in $G$. Let $v_{i_{1}}$, $v_{i_{2}}$, ..., $v_{i_{\delta(H)}}$ be neighbors of $v_{1}$ in $H$. Set $T_{j}^{'}=(v_{1}v_{i_{j}})^{u_{1}} \cup (v_{1}v_{i_{j}})^{u_{2}} \cup (v_{1}v_{i_{j}})^{u_{3}} \cup G^{v_{i_{j}}}$, $j=1$, $2$, ..., $\delta(H)$. It is clear that $T_{1}^{v_{1}}$, $T_{2}^{v_{1}}$, ..., $T_{\kappa_{3}(G)}^{v_{1}}$, $T_{1}^{'}$, $T_{2}^{'}$, ..., $T_{\delta(H)}^{'}$ are connected graphs containing $S$ and are pairwise disjoint except $S$. \qed

\medskip

From Lemma~\ref{11} to Lemma~\ref{14}, without loss of generality, we discuss all positions of three vertices of $G \square H$. Hence, theorem~\ref{15} is obvious.

~~~\\
\textbf{Theorem~\ref{15}.}~Let $G$ and $H$ be nontrivial connected graphs. Then $\kappa_{3}(G \square H) \geq min \{ \kappa_{3}(G)+ \delta(H)$, $\kappa_{3}(H)+ \delta(G)$, $\kappa(G) + \kappa(H)-1\}$.

\begin{ex}\label{18}
Let $a$ and $b$ be integers such that $a \geq 1$ and $b \geq 2$. Then $\kappa_{3}(K_{a+1} \square K_{b})=\kappa_{3}((K_{a} \vee \overline{K_{2}}) \square K_{b} )=a+b-2$.
\end{ex}

\proof
According to theorem~\ref{6}, $\kappa_{3}(K_{a+1} \square K_{b})$, $\kappa_{3}((K_{a} \vee \overline{K_{2}}) \square K_{b}) \leq a+b-2$. It is easy to see that $\kappa (K_{b})=b-1$; $\kappa_{3}(K_{b})=b-2$, if $b\geq 3$; and $\kappa(K_{a} \vee \overline{K_{2}}) = \kappa_{3}(K_{a} \vee \overline{K_{2}}) = a$. From lemma~\ref{11}-\ref{14}, we can see that $\kappa_{3}(K_{a+1} \square K_{b})$, $\kappa_{3}((K_{a} \vee \overline{K_{2}}) \square K_{b}) \geq a+b-2$.
\qed

Due to theorem~\ref{7}, $\delta(G) \geq \kappa(G) \geq \kappa_{3}(G)$, $\delta(H) \geq \kappa(H) \geq \kappa_{3}(H)$. Hence, theorem~\ref{15} induces theorem~\ref{1}.

\begin{cor}\label{17}
Let $G$ be a nontrivial connected graph, and $H$ an $l$-connected graph, where $1 \leq l \leq 5$. The following assertions hold;

(i) if $\kappa(G)=\kappa_{3}(G)$, then $\kappa_{3}(G \square H) \geq \kappa_{3}(G) + l -1$. Moreover, the bound is sharp;

(ii) if $\kappa(G) > \kappa_{3}(G)$, then $\kappa_{3}(G \square H) \geq \kappa_{3}(G) + l$. Moreover, the bound is sharp.
\end{cor}

\proof
Due to theorem~\ref{7}, $\kappa_{3}(H) \geq l-1$. Hence, $\kappa_{3}(H) + \delta(G) \geq \kappa(G) + l-1 $. And theorem~\ref{15} induces this corollary. Example~\ref{18} guarantees the sharpness.
\qed

Obviously, corollary~\ref{17} induces theorem~\ref{2} and confirms conjecture~\ref{3}. However, corollary~\ref{17} does not answer the question whether conjecture~\ref{3} still holds if $H$ is $l$-connected for $l \geq 6$. And this is what will discuss in the next section.

\section{Another lower bound}

\begin{lem}\label{19}
Let $0 \leq t \leq \lfloor \frac{l}{2} \rfloor$, and $S = \{ (u_{1}, v_{1}), (u_{1}, v_{2}), (u_{1}, v_{3})\}$. If there exists an $(l, t)$-reduced-path-bundle $\{ P_{1}, P_{2}, ..., P_{l}\} \cup \{ M_{1}, M_{2}, ..., M_{l-2t} \}$ in $H$ such that for any $1 \leq i \leq l-2t$, the terminal vertex of $M_{i}$ is on $P_{t+i}$, then

\begin{equation}
\kappa(S) \geq
   \begin{cases}
   l + \delta_{1} - 1& \text{if $\delta_{1} \geq t - 2$}ㄛ\\
   l + \delta_{1} - \lceil \frac{t - \delta_{1}}{2} \rceil& \text{if $\delta_{1} \leq t - 3$}.
   \end{cases}
\end{equation}

\end{lem}

\proof
For any $j \in \{ 1, 2, ..., t \}$, $P_{j}$ is divided into two paths by $v_{3}$, denoted by $P_{j1}$ and $P_{j2}$, where $P_{j1}$ is a $v_{1}v_{3}$-path and $P_{j2}$ is a $v_{3}v_{2}$-path. Without loss of generality, suppose $u_{1}u_{2}$, $u_{1}u_{3}$, ..., $u_{1}u_{\delta_{1}+1} \in E(G)$; $l(P_{j1})$, $l(P_{j2})$, $l(P_{l-j+1})\geq 2$, $j=1$, $2$, ..., $t-1$. Set $T_{j} = (M_{j} \cup P_{t+j})^{u_{1}}$, $j=1$, $2$, ..., $l-2t$.

Case~1~$t \leq 2$.

Set $T_{j}^{'} = (u_{1}u_{j+1})^{v_{1}} \cup (u_{1}u_{j+1})^{v_{2}} \cup (u_{1}u_{j+1})^{v_{3}} \cup H^{u_{j+1}} $, $j=1$, $2$, ..., $\delta_{1}$. If $t=0$, then $\kappa(S) \geq l-2t +\delta_{1} = l+\delta_{1}$. If $t=1$, then set $T^{''} = P_{1}^{u_{1}} $. So $\kappa(S) \geq l-2t + \delta_{1} + 1 = l+\delta_{1}-1$. If $t=2$, then set $T_{1}^{''} = (P_{l} \cup P_{11})^{u_{1}}$, $T_{2}^{''} = (P_{21} \cup P_{12})^{u_{1}}$, and $T_{3}^{''} = (P_{22} \cup P_{l-1})^{u_{1}}$. So $\kappa(S) \geq l-2t + \delta_{1} +3 =l+\delta_{1}-1$.

Case~2~$t \geq 3$.

Since $t \leq \lfloor \frac{l}{2} \rfloor$, $l \geq 6$. Suppose $v_{1}v_{i_{1}}$, $v_{1}v_{i_{2}}$, $v_{3}v_{i_{3}}$, $v_{3}v_{i_{4}}$, $v_{2}v_{i_{5}}$, $v_{2}v_{i_{6}} \in E(G)$, as shown in figure~\ref{f9}. It is easy to see that $\kappa((H - \{v_{1}, v_{2}, v_{3} \}) \cup \{ v_{i_{1}}v_{i_{2}}, v_{i_{3}}v_{i_{4}}, v_{i_{5}}v_{i_{6}}\}) \geq l-3 \geq 3$. If $\{ v_{i_{1}}v_{i_{2}}, v_{i_{3}}v_{i_{4}}, v_{i_{5}}v_{i_{6}}\}$ is not an edge cut of $(H - \{v_{1}, v_{2}, v_{3} \}) \cup \{ v_{i_{1}}v_{i_{2}}, v_{i_{3}}v_{i_{4}}, v_{i_{5}}v_{i_{6}}\}$, then according to lemma~\ref{5}, there is a cycle containing $\{ v_{i_{1}}v_{i_{2}}, v_{i_{3}}v_{i_{4}}, v_{i_{5}}v_{i_{6}}\}$ in $(H - \{v_{1}, v_{2}, v_{3} \}) \cup \{ v_{i_{1}}v_{i_{2}}, v_{i_{3}}v_{i_{4}}, v_{i_{5}}v_{i_{6}}\}$, which without loss of generality, we suppose to be $v_{i_{1}}v_{i_{2}}C_{1}v_{i_{3}}v_{i_{4}}C_{2}v_{i_{5}}v_{i_{6}}C_{3}v_{i_{1}}$, where $C_{1}$, $C_{2}$, $C_{3}$ are paths. Set $T_{1}^{'} = P_{12}^{u_{1}} \cup (u_{1}u_{2})^{v_{1}} \cup (u_{1}u_{2})^{v_{i_{3}}} \cup (v_{1}v_{i_{2}})^{u_{2}} \cup C_{1}^{u_{2}}$, $T_{2}^{'} = P_{l}^{u_{1}} \cup (u_{1}u_{2})^{v_{3}} \cup (u_{1}u_{2})^{v_{i_{5}}} \cup (v_{3}v_{i_{4}})^{u_{2}} \cup C_{2}^{u_{2}}$, and $T_{3}^{'} = P_{11}^{u_{1}} \cup (u_{1}u_{2})^{v_{2}} \cup (u_{1}u_{2})^{v_{i_{1}}} \cup (v_{2}v_{i_{6}})^{u_{2}} \cup C_{3}^{u_{2}}$. Hence, we can find $3$ internally disjoint $S$-trees in $(\underset{v \in V(P_{1} \cup P_{l})}{\cup} (u_{1}u_{2})^{v}) \cup (P_{1} \cup P_{l})^{u_{1}} \cup H^{u_{2}}$.

\begin{figure}[hptb]
  \centering
  \includegraphics[width=11cm]{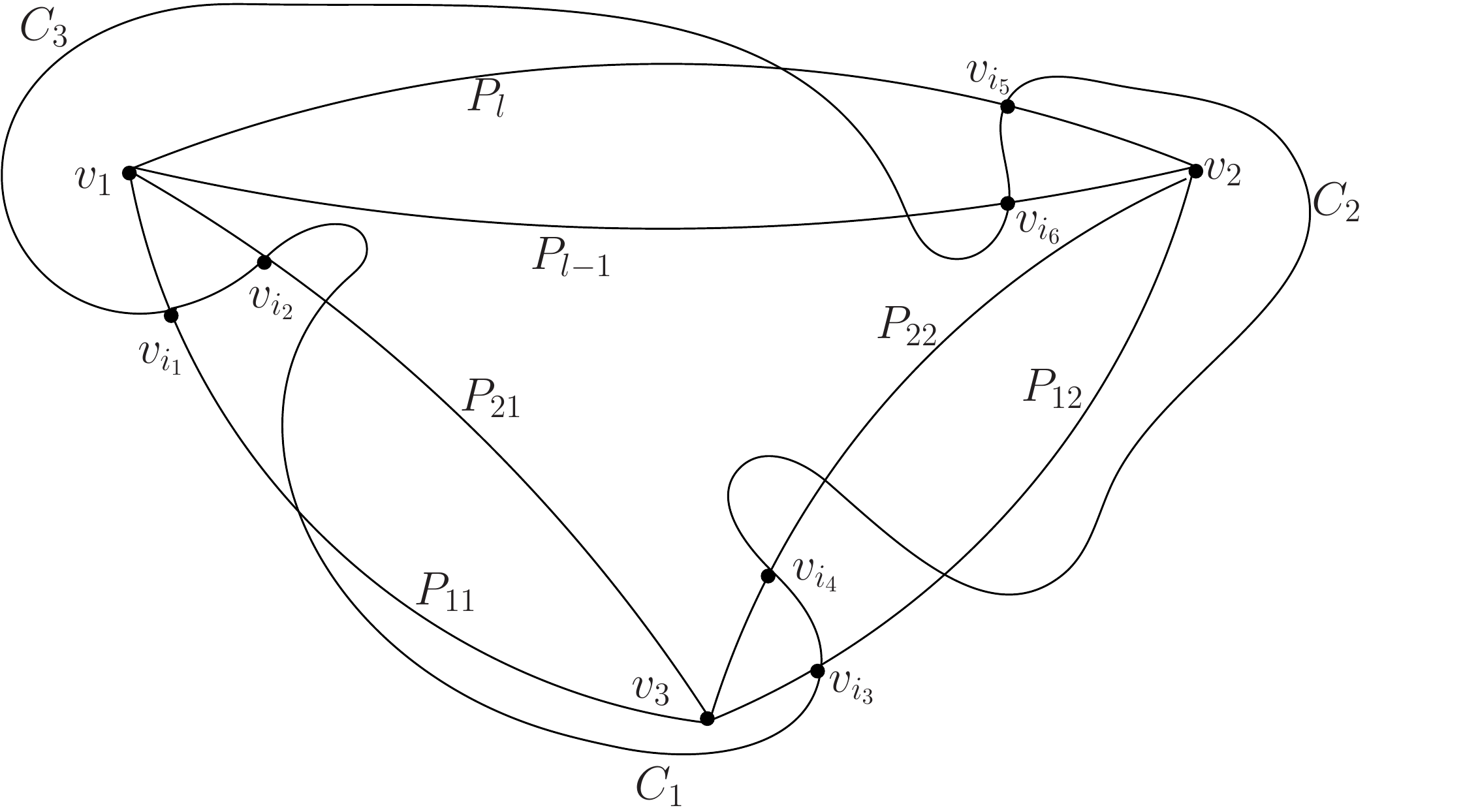}\\
  \caption{}\label{f9}
\end{figure}

Case~2.1~$l \geq 7$.

Since $\kappa((H - \{v_{1}, v_{2}, v_{3} \}) \cup \{ v_{i_{1}}v_{i_{2}}, v_{i_{3}}v_{i_{4}}, v_{i_{5}}v_{i_{6}}\}) \geq l-3 \geq 4$, $\{ v_{i_{1}}v_{i_{2}}, v_{i_{3}}v_{i_{4}}, v_{i_{5}}v_{i_{6}}\}$ cannot be an edge cut of $(H - \{v_{1}, v_{2}, v_{3} \}) \cup \{ v_{i_{1}}v_{i_{2}}, v_{i_{3}}v_{i_{4}}, v_{i_{5}}v_{i_{6}}\}$.

Similar to the discussion above, if $\delta_{1} \geq t-2$, without loss of generality, suppose we can find $3$ internally disjoint $S$-trees in $( \underset{v \in V(P_{j} \cup P_{l-j+1})}{\cup}  (u_{1}u_{1+j})^{v}) \cup (P_{j} \cup P_{l-j+1})^{u_{1}} \cup H^{u_{1+j}}$ respectively, $j=1$, $2$, ..., $t-2$. Set $T_{j}^{''} = (u_{1}u_{j})^{v_{1}} \cup (u_{1}u_{j})^{v_{2}} \cup (u_{1}u_{j})^{v_{3}} \cup H^{u_{j}}$, $j=t$, $t+1$, ..., $\delta_{1}+1$; $T_{1}^{'''} = (P_{l-t+2} \cup P_{t-1,1})^{u_{1}}$; $T_{2}^{'''} = (P_{t1} \cup P_{t-1,2})^{u_{1}}$; and $T_{3}^{'''} = (P_{t2} \cup P_{l-t+1})^{u_{1}}$. Hence, $\kappa(S) \geq l-2t + 3(t-2) + \delta_{1}-t+2 + 3 = l+\delta_{1}-1$.

If $\delta_{1} \leq t-3$, without loss of generality, suppose we can find $3$ internally disjoint $S$-trees in $( \underset{v \in V(P_{j} \cup P_{l-j+1})}{\cup}  (u_{1}u_{1+j})^{v}) \cup (P_{j} \cup P_{l-j+1})^{u_{1}} \cup H^{u_{1+j}}$ respectively, $j=1$, $2$, ..., $\delta_{1}$. Clearly, there exist at least $\lfloor \frac{3(t-\delta_{1})}{2} \rfloor$ internally disjoint $S$-trees in $(\overset{t}{\underset{j=\delta_{1}+1}{\cup}} (P_{j1} \cup P_{j2} \cup P_{l-j+1}))^{u_{1}}$. Hence, $\kappa(S) \geq l-2t + 3\delta_{1} + \lfloor \frac{3(t-\delta_{1})}{2}\rfloor= l+\delta_{1}-\lceil \frac{t-\delta_{1}}{2}\rceil$.

Case~2.2~$l=6$.

Since $3 \leq t \leq \lfloor \frac{l}{2}\rfloor$, $t=3$. If $\{ v_{i_{1}}v_{i_{2}}, v_{i_{3}}v_{i_{4}}, v_{i_{5}}v_{i_{6}}\}$ is not an edge cut of $(H - \{v_{1}, v_{2}, v_{3} \}) \cup \{ v_{i_{1}}v_{i_{2}}, v_{i_{3}}v_{i_{4}}, v_{i_{5}}v_{i_{6}}\}$, then similar to the case that $l \geq 7 $ and $\delta_{1} \geq t-2$, $\kappa(S) \geq 5 +\delta_{1}$. From now on, suppose the contrary, that is $\kappa(S) \leq 4+\delta_{1}$. Then, $u_{i_{1}}u_{i_{2}}$, $u_{i_{3}}u_{i_{4}}$, $u_{i_{5}}u_{i_{6}} \in E(H)$; and $\{ v_{i_{1}}v_{i_{2}}, v_{i_{3}}v_{i_{4}}, v_{i_{5}}v_{i_{6}}\}$ separates $H - \{v_{1}, v_{2}, v_{3} \} $ into two components $C_{1}$ and $C_{2}$, where we suppose $v_{i_{1}}, v_{i_{3}}, v_{i_{5}} \in V(C_{1})$, and $v_{i_{2}}, v_{i_{4}}, v_{i_{6}} \in V(C_{2})$. Suppose $l(P_{31}) \geq 2$, and $v_{i_{7}}$ is adjacent to $v_{1}$ in $P_{31}$. Since $\kappa(S) \leq 4+\delta_{1}$, so $v_{i_{1}}v_{i_{7}}$, $v_{i_{2}}v_{i_{7}} \in E(G)$, which induces that $ H - \{v_{1}, v_{2}, v_{3} \}- \{ v_{i_{1}}v_{i_{2}}, v_{i_{3}}v_{i_{4}}, v_{i_{5}}v_{i_{6}}\}$ is connected, a contradiction. Hence, $l(P_{31})=1$. Similarly, $l(P_{32}) = l(P_{4})=1$. Suppose $l(P_{11}) \geq 3$, and $v_{i_{7}}$ and $v_{i_{8}}$ are adjacent to $v_{3}$ in $P_{11}$ and $P_{21}$ respectively. Since $\kappa(S) \leq 4 + \delta_{1}$, $v_{i_{7}}v_{i_{8}} \in E(H)$, which induces that $ H - \{v_{1}, v_{2}, v_{3} \}- \{ v_{i_{1}}v_{i_{2}}, v_{i_{3}}v_{i_{4}}, v_{i_{5}}v_{i_{6}}\}$ is connected, a contradiction. Hence, $l(P_{11})=2$. Similarly, $l(P_{12})= l(P_{21})= l(P_{22})= l(P_{5})= l(P_{6})=2$. Suppose $v_{i_{1}}v_{2} \in E(G)$. Then in $H$, $\kappa(\{v_{1}, v_{2}, v_{3} \}) \geq 5$, and according to lemma~\ref{14}, $\kappa(S) \geq 5 + \delta_{1}$, a contradiction. Hence, $v_{i_{1}}v_{2} \notin E(G)$. Since $d(v_{i_{1}}) \geq 6$, there exists a vertex $v_{i_{9}} \in V(C_{1})-\{v_{i_{1}}, v_{i_{3}}, v_{i_{5}} \}$ such that $v_{i_{9}}v_{i_{1}} \in E(H)$. Since $\kappa(H)=6$, so there exists a $6$-fan in $H$ from $v_{i_{9}}$ to $\{ v_{1}, v_{2}, v_{3}, v_{i_{1}}, v_{i_{3}}, v_{i_{5}}\}$, disjoint from $C_{2}$. Then in $H$, $\kappa(\{ v_{1}, v_{2}, v_{3}\}) \geq 5$, also a contradiction. \qed

\begin{pro}\label{20}

Let $S= \{ (u_{1}, v_{1}), (u_{1}, v_{2}), (u_{1}, v_{3})\}$. Then

\begin{equation}
\kappa(S) \geq
   \begin{cases}
   l + \delta_{1} - 1& \text{if $\delta_{1} \geq \lfloor \frac{l}{2} \rfloor - 2$}, \\
   l + \delta_{1} - \lceil \frac{\lfloor \frac{l}{2} \rfloor - \delta_{1}}{2} \rceil& \text{if $\delta_{1} \leq \lfloor \frac{l}{2} \rfloor - 3$}.
   \end{cases}
\end{equation}

\end{pro}

\proof
Since $H$ is $l$-connected, so according to theorem~\ref{7}, for some $0 \leq t \leq \lfloor \frac{l}{2} \rfloor $, there exists an $(l, t)$-reduced-path-bundle $\{ P_{1}, P_{2}, ..., P_{l}\} \cup \{ M_{1}, M_{2}, ..., M_{l-2t} \}$ such that for any $1 \leq i \leq l-2t$, the terminal vertex of $M_{i}$ is on $P_{t+i}$. When $\delta_{1} \geq \lfloor \frac{l}{2}\rfloor - 2$, Since $t \leq \lfloor \frac{l}{2} \rfloor$, $\delta_{1} \geq t-2$. Hence, due to lemma~\ref{19}, $\kappa(S) \geq l+\delta_{1}-1$. When $\delta_{1} \leq \lfloor \frac{l}{2}\rfloor -3$, if $\delta_{1} \geq t-2$, then $\kappa(S) \geq l+\delta_{1}-1 \geq l+\delta_{1}- \lceil \frac{\lfloor \frac{l}{2} \rfloor - \delta_{1}}{2} \rceil$; if $\delta_{1} \leq t-3$, according to lemma~\ref{19}, $\kappa(S) \geq l + \delta_{1} - \lceil \frac{ t - \delta_{1}}{2} \rceil \geq l + \delta_{1} - \lceil \frac{\lfloor \frac{l}{2} \rfloor - \delta_{1}}{2} \rceil$. \qed

~\\
\textbf{Theorem~\ref{16}}~Let $G$ be a nontrivial connected graph, and let $H$ be an $l$-connected graph. The following assertions hold:

(i) if $\kappa(G)=\kappa_{3}(G)$ and $1 \leq l \leq 7$, then $\kappa_{3}(G \square H) \geq \kappa_{3}(G) + l -1$. Moreover, the bound is sharp;

(ii) if $\kappa(G) > \kappa_{3}(G)$ and $1 \leq l \leq 9$, then $\kappa_{3}(G \square H) \geq \kappa_{3}(G) + l$. Moreover, the bound is sharp.

\proof
Let $S$ be a $3$-vertex-set of $V(G \square H)$. If there does not exist a vertex $u \in V(G)$ such that $S \subseteq V(H^{u})$, then according to lemma~\ref{11}-lemma~\ref{14}, $\kappa(S) \geq \kappa_{3}(G) + \delta(H) \geq \kappa_{3}(G) + l$ or $\kappa(S) \geq \kappa(G)+l-1$. Next, without loss of generality, suppose $S= \{(u_{1}, v_{1}), (u_{1}, v_{2}), (u_{1}, v_{3}) \}$. If $l \leq 7$, then $\delta_{1} \geq \lfloor \frac{l}{2} \rfloor -2$, which according to proposition~\ref{20}, induces that $\kappa(S) \geq l+\delta_{1}-1$. If $l=8$ or $9$ and $\kappa(G) > \kappa_{3}(G)$, then $\delta_{1} \geq 2 = \lfloor \frac{l}{2} \rfloor -2 $, which according to proposition~\ref{20}, induces that $\kappa(S) \geq l+\delta_{1}-1 \geq \kappa_{3}(G) +l$. Example~\ref{18} guarantees the sharpness. \qed

An $S$-tree $T$ in $G$ is said to be minimal if for any $S$-tree $T^{'}$, $V(T^{'}) \subseteq V(T)$ and $\partial{(S)} \cap E(T^{'}) \subseteq \partial{(S)} \cap E(T)$ induce $V(T^{'}) = V(T)$ and $\partial{(S)} \cap E(T^{'}) = \partial{(S)} \cap E(T)$, where $\partial(S)$ denotes the set of edges with one vertex in $S$.

\begin{ex}\label{21}
Let $a$, $b$, $c$ be integers such that $1 \leq a \leq b \leq c$. Then,

\begin{equation}
\kappa_{3}(K_{a,b,c}) =
   \begin{cases}
   a+b& \text{if $a+b \leq c$}, \\
   \lfloor \frac{a+b+c}{2}\rfloor& \text{if $a+b \geq c+1$}.
   \end{cases}
\end{equation}

\end{ex}

\proof
It is easy to see that $\kappa_{3}(K_{1,1,1})=1$ and $\kappa_{3}(K_{1,1,c})=2$ for $c \geq 2$. Next, suppose $a+b \geq 3$. Let $V_{1}=\{ v_{11}, v_{21}, ..., v_{a1}\}$, $V_{2}=\{ v_{12}, v_{22}, ..., v_{b2}\}$ and $V_{3}=\{ v_{13}, v_{23}, ..., v_{c3}\}$ be the maximal independent sets of $K_{a,b,c}$. Let $S$ be a $3$-vertex-set of $V(K_{a,b,c})$. If $S \subseteq V_{1} \cup V_{2}$, then $\kappa(S) \geq c+\kappa_{3}(K_{a,b})$; If $S \subseteq V_{1} \cup V_{3}$, then $\kappa(S) \geq b+\kappa_{3}(K_{a,c})$; If $S \subseteq V_{2} \cup V_{3}$, then $\kappa(S) \geq a+\kappa_{3}(K_{b,c})$. Moreover, equalities can hold in the three inequalities above. Next, without loss of generality, suppose $S = \{ v_{11}, v_{12}, v_{13}\}$. Then figure~\ref{f10} shows three types of minimal $S$-trees, namely A, B and C. Let $\mathscr{T} $ be the set of $\kappa(S)$ internally disjoint minimal $S$-trees, containing $u$ trees of type A (or A trees), $v$ trees of type B (or B trees) and $w$ trees of type C (or C trees). We notice that $v+2w \leq 3$ and $2u+v \leq a+b+c-3$. Hence, $\kappa(S)=u+v+w \leq \lfloor \frac{a+b+c}{2} \rfloor$. Besides, $\kappa(S) \leq d(v_{13}) = a+b$. Set $T_{1} = (v_{11}v_{12}) \cup (v_{12}v_{13})$ and $T_{2} = (v_{11}v_{13})\cup(v_{11}v_{23})\cup(v_{12}v_{23})$. If $a+b \geq c+1$, then it is easy to see that there exist $\lfloor \frac{a-1+b-1+c-2}{2}\rfloor = \lfloor \frac{a+b+c}{2}\rfloor - 2$ internally disjoint A trees in $K_{a,b,c}-v_{23}-\{ v_{11}v_{12}, v_{11}v_{13}, v_{12}v_{13}\}$. So $\kappa(S)=\lfloor \frac{a+b+c}{2}\rfloor $. If $a+b \leq c$, then it is easy to see that there exist $a-1+b-1=a+b-2$ internally disjoint A trees in $K_{a,b,c}-v_{23}-\{ v_{11}v_{12}, v_{11}v_{13}, v_{12}v_{13}\}$. So $\kappa(S)=a+b$. Finally, because of lemma~\ref{9}, the problem is solved by comparing sizes. \qed

\begin{figure}[hptb]
  \centering
  \includegraphics[width=11cm]{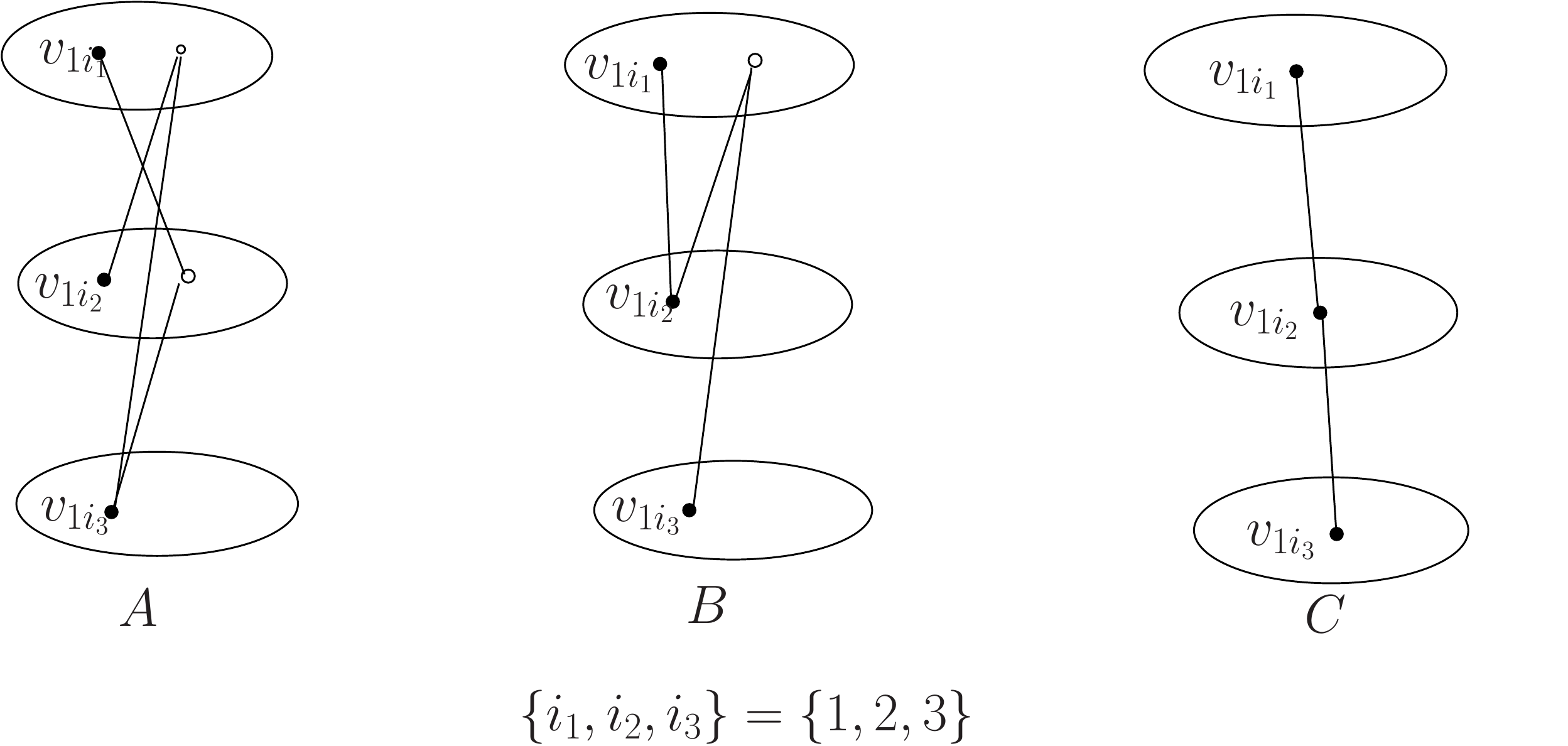}\\
  \caption{}\label{f10}
\end{figure}

\begin{figure}[htpb]
  \centering
  \includegraphics[width=10cm]{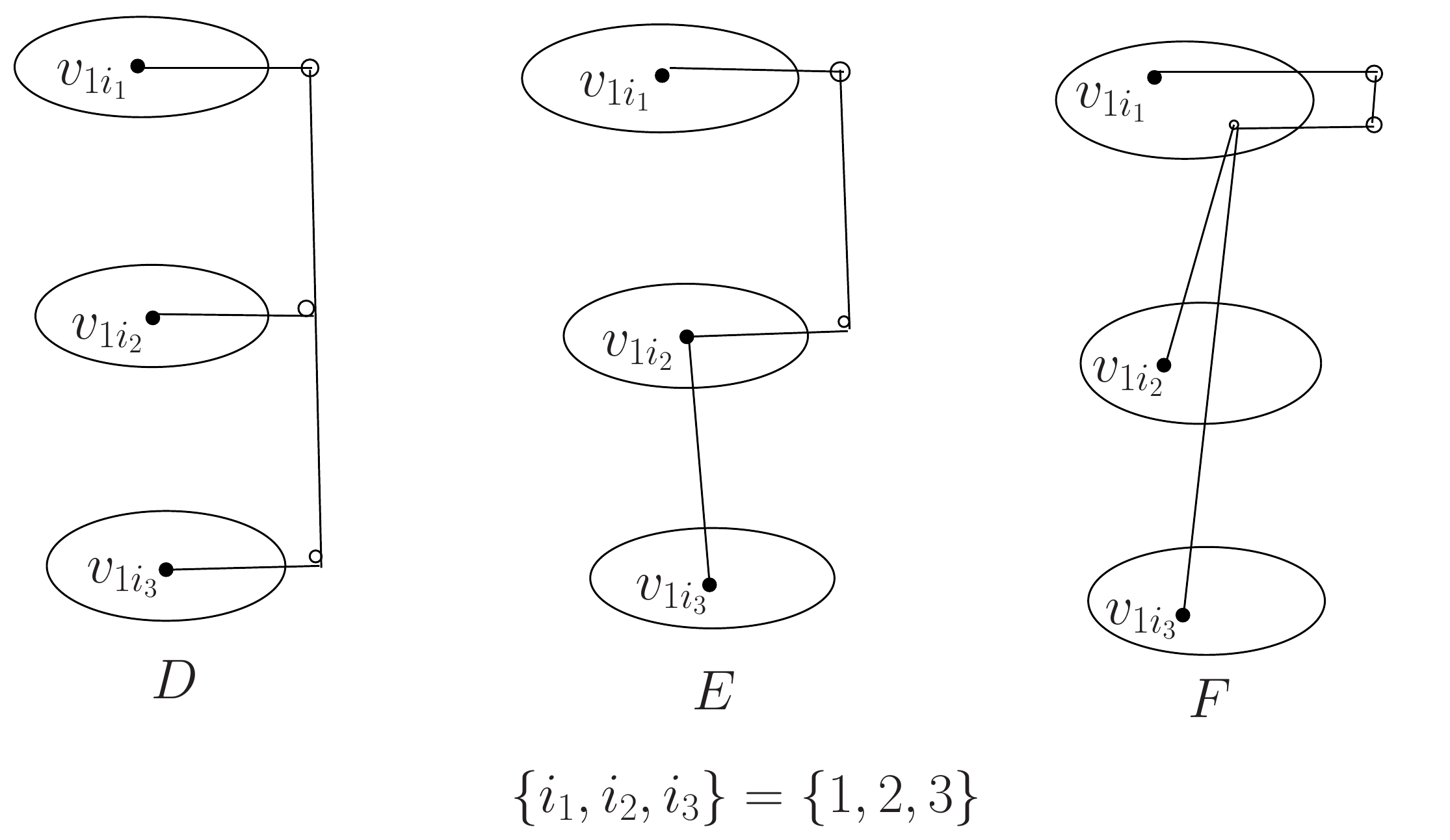}\\
  \caption{}\label{f12}
\end{figure}
\begin{ex}\label{22}
Let $a$ and $b$ be integers such that $a\geq 1$ and $b \geq 2$. Then,

\begin{equation}
\kappa_{3}(K_{b} \square K_{a,a,a}) =
   \begin{cases}
   2a + b -2& \text{if $b \geq a-1$}, \\
   \lfloor \frac{3a+3b-3}{2}\rfloor& \text{if $b \leq a-2$}.
   \end{cases}
\end{equation}

and

\begin{equation}
\kappa_{3}(K_{b} \square K_{a,a+1,a+1}) =
   \begin{cases}
   2a + b -1& \text{if $b \geq a-1$}, \\
   \lfloor \frac{3a+3b-1}{2}\rfloor& \text{if $b \leq a-2$}.
   \end{cases}
\end{equation}

\end{ex}

\proof
Since the proofs are quite similar for $K_{b} \square K_{a,a,a}$ and $K_{b} \square K_{a,a+1,a+1}$, we give our proof only for the latter. Let $V(K_{b}) = \{ u_{1}, u_{2}, ..., u_{b}\}$; Let $V_{1}=\{ v_{11}, v_{21}, ..., v_{a1}\}$, $V_{2}=\{ v_{12}, v_{22}, ..., v_{a+1, 2}\}$ and $V_{3}=\{ v_{13}, v_{23}, ..., v_{a+1, 3}\}$ be the maximal independent sets of $K_{a,a+1,a+1}$; and Let $S$ be a $3$-vertex-set of $V(K_{b} \square K_{a,a+1,a+1})$. Because of theorem~\ref{6}, $\kappa_{3}(K_{b} \square K_{a,a+1,a+1}) \leq 2a+b-1$. If for any $u \in V(G)$, $S \nsubseteq V(H^{u})$, then due to lemma~\ref{11}-\ref{14}, $\kappa(S) \geq 2a+b-1$; If there exists a vertex $u \in V(G)$ such that for some $i \in \{ 1, 2, 3\}$, $S \subseteq (V(H)-V_{i})^{u}$, then according to the proof of lemma~\ref{14}, $\kappa(S) \geq 2a+b-1$. Next, without loss of generality, let $S = \{ (u_{1}, v_{11}), (u_{1}, v_{12}), (u_{1}, v_{13})\}$. Then according to proposition~\ref{20}, we have

\begin{equation}
\kappa(S) \geq
   \begin{cases}
   2a + b -1& \text{if $b \geq a-1$}, \\
   \lfloor \frac{3a+3b-1}{2}\rfloor& \text{if $b \leq a-2$}.
   \end{cases}
\end{equation}

Hence, when $b \geq a-1$, $\kappa_{3}(K_{b} \square K_{a,a+1,a+1}) = 2a+b-1$. Next, we claim that when $b \leq a-2$, $\kappa(S) \leq \lfloor \frac{3a+3b-1}{2}\rfloor$.

Let $F$ be the graph obtained from $K_{b} \square K_{a,a+1,a+1}$ by joining every pair of nonadjacent vertices in $\{(u,v) | u \in V(K_{b})-u_{1}, v \in V(K_{a,a+1,a+1}) \}$. Obviously, the maximum number of internally disjoint $S$-trees in $F$ is not less than the maximum number of internally disjoint $S$-trees in $K_{b} \square K_{a,a+1,a+1}$. In $F$, all types of minimal $S$-trees are shown in figure~\ref{f10} and figure~\ref{f12}. Let $\mathscr{T}$ be the set with maximum number of internally disjoint minimal $S$-trees in F and as many F trees as possible. Let the number of A, B, C, D, E and F trees $\mathscr{T}$ contains be $u$, $v$, $w$, $x$, $y$ and $z$, denoted by $\mathscr{T} = uA+vB+wC+xD+yE+zF$. It is easy to see that one A tree and one D tree can be replaced by two F trees, denoted by $A+D \rightarrow 2F$. So A trees and D trees can not be both in $\mathscr{T}$. Similarly, we have $A+E \rightarrow B+F$. If $x >0$ or $y>0$, then $u=0$. We notice that $3x+2y+z \leq 3(b-1)$ and $v+2w+y \leq 3$, so $|\mathscr{T}|=v+w+x+y+z \leq 3b < \lfloor \frac{3a+3b-1}{2}\rfloor$. If $x=y=0$, then $z=3(b-1)$ and $|\mathscr{T}| = 3(b-1) + \kappa_{3}(K_{a-b+1,a-b+2,a-b+2}) = \lfloor \frac{3a+3b-1}{2}\rfloor$. \qed

Also, we can similarly prove that when $a \geq 4$, $\kappa_{3}(P_{3} \square K_{a,a,a}) = \kappa_{3}(K_{2} \square K_{a,a,a}) = \lfloor \frac{3a + 3}{2} \rfloor$, and $\kappa_{3}(P_{3} \square K_{a,a+1,a+1}) = \kappa_{3}(K_{2} \square K_{a,a+1,a+1}) = \lfloor \frac{3a + 5}{2} \rfloor$. In theorem~\ref{16}, let $G=P_{3}$, and $H=K_{a,a,a}$ or $K_{a,a+1,a+1}$. Then $l \geq 8$, and $\kappa_{3}(G \square H) < \kappa_{3}(G) + l-1$. When $a \geq 5$, $\kappa_{3}(K_{3} \square K_{a,a,a}) =  \lfloor \frac{3a + 6}{2} \rfloor$, and $\kappa_{3}(K_{3} \square K_{a,a+1,a+1}) =  \lfloor \frac{3a + 8}{2} \rfloor$. In theorem~\ref{16}, let $G=K_{3}$, and $H=K_{a,a,a}$ or $K_{a,a+1,a+1}$. Then $l \geq 10$, and $\kappa_{3}(G \square H) < \kappa_{3}(G) + l$.

\end{document}